\documentclass[10pt]{amsart}
\usepackage{latexsym}
\usepackage{graphics}
\usepackage{amssymb}
\usepackage{amsbsy}
\usepackage{amsmath}
\newtheorem{Lemma}{Lemma} \par
\newtheorem{Theorem}[Lemma]{Theorem} \par
\newtheorem{Remark}[Lemma]{Remark} \par
\newcommand{\RKA}{{\mathcal{A}}}
\newcommand{\R}{{\mathbb R}}
\newcommand{\N}{{\mathbb N}}
\newcommand{\multsp}{\,}
\begin{document}
\title{Runge-Kutta Methods, Trees, and {\em Mathematica}}
\author{Folkmar Bornemann}
\address{Munich University of
Technology, 80290 Munich, Germany}
\email{bornemann@ma.tum.de}
\date{November 4, 2002}
\keywords{numerical solution of ODEs, Runge-Kutta methods, recursive
representation of rooted trees,
Butcher's theorem, automatic generation of order conditions,
computer algebra systems.}
\subjclass{65-01,65L06,65Y99}
\begin{abstract}
This paper presents a simple and elementary proof of But\-cher's
theorem on the order conditions of Runge-Kutta methods. It is
based on a recursive definition of rooted trees and avoids
combinatorial tools such as labelings and Fa\`{a} di Bruno's formula.
This strictly recursive approach can easily and elegantly be
implemented using modern computer algebra systems like {\em Mathematica} for
automatically generating the order conditions. The full, but short
source code is presented and applied to some instructive examples.
\end{abstract}
\maketitle
\pagestyle{myheadings}
\thispagestyle{plain}
\markboth{F. BORNEMANN}{RUNGE-KUTTA METHODS, TREES, AND {\em MATHEMATICA}}
\section{Introduction}

A first step towards the construction of Runge-Kutta
methods is the calculation of the order conditions that the
coefficients have to obey. In the old days they were obtained
by expanding the error term in a Taylor series by hand, a procedure
which for higher orders sooner or later runs into difficulties
because of the largely increasing combinatorial complexity. It
was a  major break-through when Butcher \cite{Butc63} published 1963 his
result of systematically describing order conditions by rooted trees.
The proof of this result has evolved very much in meantime,
mainly under the influence of Butcher's later work \cite{Butc87}
and the contributions of Hairer and Wanner \cite{HaWa74,HNWa87}.
In this paper we will present a simple and elementary proof of
Butcher's theorem by using very consequently the recursive
structure of rooted trees. This way we avoid lengthy calculations
of combinatorial coefficients, the use of tree-labelings, or
Fa\`{a} di Bruno's formula as in \cite{Butc87,HNWa87}. Our
proof is very similar in spirit to the presentation of $B$-series by
Hairer in Chapter 2 of his lecture notes \cite{Hair99}.

As early as 1976 Jenks \cite{Jenk76} posed the problem of
automatically generating order conditions for Runge-Kutta methods
using computer algebra systems, but no replies were received.
In 1988 Keiper of Wolfram Research concluded that the
method of automatically calculating Taylor expansions by brute
force was bound to be very inefficient. Naturally he turned to
the elegant results of Butcher's. Utilizing them, he wrote the
{\em Mathematica} package {\tt Butcher.m}, which has been available as
part of the standard distribution of {\em Mathematica} since then.
This package was later considerably improved by Sofroniou
\cite{Sofr94} and offers a lot of sophisticated tools.

While teaching the simple proof of Butcher's result in a first
course on numerical ODEs, the author realized that the underlying
recursive structure could also be exploited for a simple and elegant
computer implementation. This approach differs from the work of Sofroniou
in various respects. We will present the full source code in {\em Mathematica}
and some applications.

\section{Runge-Kutta methods}

The Runge-Kutta methods are one-step discretizations of
initial-value problems for systems of $d$ ordinary differential
equations,
\[
x' = f(t,x),\qquad x(t_0) = x_0,
\]
where the right-hand side $f: [t_0,T]\times \Omega \subset \R \times \R^d
\to\R^d$ is assumed to be sufficiently smooth. The {\em continuous evolution} $x(t) = \Phi^{t,t_0} x_0$ of the initial-value
problem is
approximated in steps of length $\tau$ by $\Psi^{t+\tau,t} x \approx
\Phi^{t+\tau,t}x$. This {\em discrete evolution} $\Psi^{t+\tau,t}x$ is defined as an approximation
of the integral-equation representation
\[
\Phi^{t+\tau,t} x = x + \int_t^{t+\tau} f(\sigma,\Phi^{\sigma,t}
x) \,d\sigma
\]
by appropriate quadrature formulas:
\begin{equation}\label{eq:rk}
\begin{array}{rcl}
\Psi^{t+\tau,t}x = x + \tau \displaystyle\sum_{i=1}^s b_i
k_i,\qquad
k_i = f\left(t+c_i\tau, x + \tau\displaystyle\sum_{j=1}^s
a_{ij} k_j\right).
\end{array}
\end{equation}
The vectors $k_i \in \R^d$, $i=1,\ldots,s$, are called {\em stages}, $s$ is
the stage number. Following the standard notation, we collect the coefficients
of the method into a matrix
and two vectors
\[
\RKA = (a_{ij})_{ij} \in \R^{s \times s},\; b = (b_1,\ldots,b_s)^T \in
\R^s,\;
c = (c_1,\ldots,c_s)^T  \in \R^s.
\]
The method is {\em explicit}, if $\RKA$ is strictly lower
triangular. The method has {\em order} $p \in \N$, if the  error term
expands to
\[
\Phi^{t+\tau,t}x - \Psi^{t+\tau,t}x = O(\tau^{p+1}).
\]
In terms of the Taylor
expansion of the error $\Phi - \Psi$, the vanishing of all lower order terms
in $\tau$ just defines the conditions which have to be satisfied
by the coefficients $\RKA$, $b$ and $c$ of a Runge-Kutta method.

If we choose $c_i = \sum_{j=1}^s a_{ij}$,
it can be shown \cite{HNWa87}, that there is no loss of generality
in considering {\em autonomous} systems only, i.e., those with
no dependence of $f$ on $t$. Doing so, the expressions
$\Phi^{t+\tau,t}x$ and $\Psi^{t+\tau,t} x$ are likewise
independent of $t$. We will write $\Phi^\tau x$ and
$\Psi^\tau x$ for short, calling them the {\em flow} and the {\em discrete flow}
of the continuous resp. the discrete system.

\section{Elementary differentials and rooted trees}
The Taylor expansions of both the phase flow $\Phi^\tau x$ and the discrete flow
$\Psi^\tau x$ are linear combinations of {\em elementary} differentials like
\[
f'''(f'f,f'f,f) = \sum_{ijklm} \frac{\partial^3 f}{\partial x_i
\partial x_j\partial x_k} \cdot \frac{\partial f_i}{\partial x_l
} f_l \cdot \frac{\partial f_j}{\partial x_m } f_m \cdot f_k .
\]
We will use the short multilinear
notation of the left hand side for the rest of the paper.

An elementary differential can be expressed uniquely by the
structure of how the subterms enter the multilinear maps. For instance,
looking at the
expression $f'''(f'f,f'f,f)$ we observe that $f'''$ must be a third
derivative, since {\em three} arguments make it {\em three}-linear. This
structure can be expressed in general by {\em rooted
trees}, e.g., \setlength{\unitlength}{1cm}
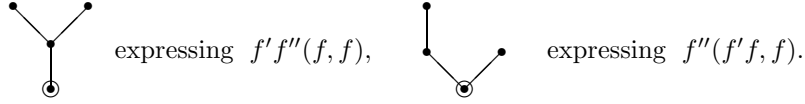
\[
\begin{picture}(5.5,2)
\put(0.75,0.5){\line(0,1){0.6}} \put(0.75,1.1){\line(1,1){0.5}}
\put(0.75,1.1){\line(-1,1){0.5}} \put(0.75,0.5){\circle{0.2}}
\put(0.75,0.5){\circle*{0.1}} \put(0.75,1.1){\circle*{0.1}}
\put(0.25,1.6){\circle*{0.1}} \put(1.25,1.6){\circle*{0.1}}
\put(3.25,1){\makebox(0,0){ expressing $\; f'f''(f,f),$}}
\end{picture}
\begin{picture}(5.5,2)
\put(0.75,0.5){\line(1,1){0.5}} \put(0.75,0.5){\line(-1,1){0.5}}
\put(0.25,1.0){\line(0,1){0.6}} \put(0.75,0.5){\circle{0.2}}
\put(0.75,0.5){\circle*{0.1}} \put(0.25,1){\circle*{0.1}}
\put(1.25,1){\circle*{0.1}} \put(0.25,1.6){\circle*{0.1}}
\put(3.5,1){\makebox(0,0){ expressing $\; f''(f'f,f).$}}
\end{picture}
\]
Every node with $n$ children denotes a $n$th derivative of $f$,
 which is applied as a multilinear
map to further elementary differentials, according to the structure
of the tree. We start reading off this structure by looking at the
root. This defines a {\em recursive} procedure, if we
observe the following: Having removed the root and its
edges, a rooted tree $\beta$ decomposes into rooted subtrees
$\beta_1,\ldots,\beta_n$ with strictly less nodes. The roots of the subtrees
$\beta_1,\ldots,\beta_n$ are exactly the $n$ children of
$\beta$'s root. This way a rooted tree $\beta$ can be
defined as the {\em unordered} list of its successors
\begin{equation}\label{eq:treedef}
\beta=[\beta_1,\ldots,\beta_n],\qquad \#\beta = 1 + \#\beta_1 +
\ldots + \#\beta_n.
\end{equation}
Here, we denote by $\#\beta$ the {\em order} of a rooted tree $\beta$, i.e.,
the number of its nodes. The root itself can be
identified with the {\em empty} list, $\odot = [\,]$.

An application of this procedure shows for the examples above that
$f'(f''(f,f))$ is expressed by $[[\odot,\odot]]$ and $f''(f'(f),f)$
is expressed by $[[\odot],\odot]$.
The reader will observe the perfect matching of parentheses and
commas. In general the relation between a rooted tree
$\beta=[\beta_1,\ldots,\beta_n]$ and its corresponding
elementary differential $f^{(\beta)}(x)$ is recursively defined
by
\begin{equation}\label{eq:diffdef}
f^{(\beta)}(x)=f^{(n)}(x) \cdot
\left(f^{(\beta_1)}(x),\ldots,f^{(\beta_n)}(x)\right).
\end{equation}
The dot of multiplication denotes the multilinear application
of the derivative to the $n$ given arguments. Due to the
symmetry of the $n$-linear map
$f^{(n)}$, the order of the subtrees
$\beta_1,\ldots,\beta_n$ does not matter, which means,
that $f^{(\beta)}$ depends in a well-defined way  on
$\beta$ as an unordered list only.

From $\odot = [\,]$ we deduce
$f^{(\odot)} = f$.
Analogously, each of the recursive definitions in the following will
have a well-defined meaning if applied to the single root $\odot =
[\,]$, mostly by using the reasonable convention that empty
products evaluate to one and empty sums to zero---a convention
that is also observed by most computer algebra systems.

\section{A simple proof of Butcher's theorem}

We are now in a position to calculate and denote the Taylor expansion of
the continuous flow $\Phi^\tau$ in a clear and compact fashion.

\begin{Lemma}
Given $f\in C^p(\Omega,\R^d)$ the flow $\Phi^\tau x$ expands to
\begin{equation}\label{eq:phitay}
\Phi^\tau x = x +  \sum_{\#\beta\leq p}
\frac{\tau^{\#\beta}}{\beta!} \alpha_\beta\, f^{(\beta)}(x) +
O(\tau^{p+1}).
\end{equation}
The coefficients $\beta!$ and $\alpha_\beta$ belonging to a rooted tree
 $\beta=[\beta_1,\ldots,\beta_n]$ are recursively defined by
\begin{equation}\label{eq:factorialdef}
\beta ! = (\#\beta)\, {\beta_1}!\cdot\ldots\cdot{\beta_n}!, \qquad
\alpha_\beta = \frac{\delta_\beta}{n !}\,
\alpha_{\beta_1}\cdot\ldots\cdot \alpha_{\beta_n}.
\end{equation}
By $\delta_\beta$ we denote the number of different ordered
tuples $(\beta_1,\ldots,\beta_n)$ which correspond to the
same unordered list $\beta =[\beta_1,\ldots,\beta_n]$.
\end{Lemma}
\begin{proof}
The assertion is obviously true for $p=0$. We proceed by
induction on $p$. Using the assertion for $p$,
the multivariate Taylor formula and the multilinearity of the
derivatives we obtain
\[
\begin{split}
f(\Phi^\tau x) &= f\left(x +  \sum_{\#\beta\leq p}
\frac{\tau^{\#\beta}}{\beta!} \alpha_\beta f^{(\beta)} +
O(\tau^{p+1})\right)\\ &= \sum_{n = 0}^p \frac{1}{n!}
f^{(n)} \!\cdot\!\! \left( \sum_{\#\beta_1\leq p}
\frac{\tau^{\#\beta_1}}{\beta_1!} \alpha_{\beta_1}
f^{(\beta_1)},\ldots,\!\sum_{\#\beta_n\leq p}
\frac{\tau^{\#\beta_n}}{\beta_n!}
 \alpha_{\beta_n}
f^{(\beta_n)} \right) + O(\tau^{p+1})\\
\end{split}
\]
\[
\begin{split}
\phantom{f(\Phi^\tau x)}&=\sum_{n=0}^p
\frac{1}{n!}\sum_{\#\beta_1+\ldots+\#\beta_n\leq p}
\frac{\tau^{\#\beta_1+\ldots+\#\beta_n}}{\beta_1! \cdot
\ldots\cdot \beta_n!} \cdot
\alpha_{\beta_1}\cdot\ldots\cdot\alpha_{\beta_n}\cdot
\\ &\qquad f^{(n)} \cdot \left( f^{(\beta_1)},\ldots,f^{(\beta_n)}\right) +
O(\tau^{p+1})\\
&= \sum_{n=0}^p \sum_{\substack{
\beta =[\beta_1,\ldots,\beta_n]\\
\#\beta \leq p+1
}}
\frac{\#\beta\cdot
\tau^{\#\beta-1}}{\beta!}\! \cdot\!
\underbrace{\frac{\delta_\beta}{n!}\,
\alpha_{\beta_1}\cdot\ldots\cdot\alpha_{\beta_n}}_{=\alpha_\beta}
\,f^{(\beta)} + O(\tau^{p+1})\\ &=
\sum_{\#\beta \leq p+1} \frac{\#\beta\cdot
\tau^{\#\beta-1}}{\beta!} \alpha_\beta\,f^{(\beta)}
+O(\tau^{p+1}).
\end{split}
\]
Plugging this into the integral form of the initial value problem
we obtain
\[
\Phi^\tau x = x + \int_0^\tau f(\Phi^\sigma x) \,d\sigma = x + \sum_{\#\beta \leq p+1}
\frac{\tau^{\#\beta}}{\beta!} \alpha_\beta\,f^{(\beta)}
+O(\tau^{p+2}),
\]
which proves the assertion for $p+1$.
\end{proof}

A likewise clear and compact expression can be calculated
for the Taylor expansion of the discrete flow.

\begin{Lemma}
Given $f\in C^p(\Omega,\R^d)$ the discrete flow $\Psi^\tau x$
expands to
\begin{equation}\label{eq:psitay}
\Psi^\tau x = x +  \sum_{\#\beta\leq p} \tau^{\#\beta}
\alpha_\beta\cdot b^T\RKA^{(\beta)}\, f^{(\beta)}(x) +
O(\tau^{p+1}).
\end{equation}
The vector $\RKA^{(\beta)} \in
\R^s$,
$\beta=[\beta_1, \ldots,\beta_n]$, is recursively defined by
\begin{equation}\label{eq:abeta}
\RKA^{(\beta)}_i = \left(\RKA\cdot \RKA^{(\beta_1)}\right)_i
\cdot\ldots\cdot \left(\RKA\cdot\RKA^{(\beta_n)}\right)_i,\qquad
i=1,\ldots,s.
\end{equation}
\end{Lemma}
\begin{proof}
Because of the definition (\ref{eq:rk}) of the discrete flow
we have to prove that the stages
$k_i$ expand to
\[
k_i = \sum_{\#\beta\leq p} \tau^{\#\beta-1}
\alpha_\beta\,\RKA^{(\beta)}_i f^{(\beta)} + O(\tau^{p}).
\]
This is obviously the case for $p=0$. We proceed by induction on
$p$. Using the assertion for $p$, the definition of the stages $k_i$,
the multivariate Taylor formula and the multilinearity of
the derivatives we obtain
\[
\begin{split}
k_i &=
f\left(x + \tau \left(\sum_{\#\beta\leq p} \tau^{\#\beta-1}
\alpha_\beta\, \left(\RKA\cdot \RKA^{(\beta)}\right)_i\,
f^{(\beta)} + O(\tau^{p}) \right)\right)\\ &= \sum_{n = 0}^p
\frac{1}{n!} f^{(n)} \cdot \left( \sum_{\#\beta_1\leq p}
\tau^{\#\beta_1} \alpha_{\beta_1}\,\left(\RKA\cdot
\RKA^{(\beta_1)}\right)_i\, f^{(\beta_1)},\ldots\right.\\
&\qquad\qquad\ldots\left., \sum_{\#\beta_n\leq p}
\tau^{\#\beta_n} \alpha_{\beta_n}\,\left(\RKA\cdot
\RKA^{(\beta_n)}\right)_i \, f^{(\beta_n)} \right) \;+
O(\tau^{p+1})\\
&=\sum_{n=0}^p
\frac{1}{n!}\sum_{\#\beta_1+\ldots+\#\beta_n\leq p}
\tau^{\#\beta_1+\ldots+\#\beta_n} \cdot
\alpha_{\beta_1}\cdot\ldots\cdot\alpha_{\beta_n}\cdot\left(\RKA\cdot \RKA^{(\beta_1)}\right)_i \cdot \\ &\;
\ldots\cdot
\left(\RKA\cdot\RKA^{(\beta_n)}\right)_i \, f^{(n)} \cdot \left(
f^{(\beta_1)},\ldots,f^{(\beta_n)}\right) + O(\tau^{p+1})
\end{split}
\]
\[
\begin{split}
&= \sum_{n=0}^p \;\sum_{\substack{
\beta =[\beta_1,\ldots,\beta_n]\\
\#\beta \leq p+1
}}
 \tau^{\#\beta-1} \cdot \underbrace{\frac{\delta_\beta}{n!}\,
\alpha_{\beta_1}\cdot\ldots\cdot\alpha_{\beta_n}}_{=\alpha_\beta}
\cdot\,\RKA^{(\beta)}_i \,f^{(\beta)} + O(\tau^{p+1})\\
&= \sum_{\#\beta \leq p+1} \tau^{\#\beta-1}
\alpha_\beta\cdot\RKA^{(\beta)}_i\,f^{(\beta)} +O(\tau^{p+1}),
\end{split}
\]
which proves the assertion for $p+1$.
\end{proof}

Comparing the coefficients of the elementary differentials in the expansions
of both the phase flow and the discrete flow, we immediately obtain Butcher's
theorem~\cite{Butc63}.

\begin{Theorem}[Butcher 1963]\label{satz:butcher}
A Runge-Kutta-method $(b,\RKA)$ is of order $p \in \N$ for all $f \in C^p(\Omega,\R^d)$,
if the order conditions
\begin{equation}\label{eq:cond}
b^T \RKA^{(\beta)} = 1/\beta !
\end{equation}
are satisfied for all rooted trees $\beta$ of order $\#\beta\leq p$.
\end{Theorem}

\section{Generating order conditions with {\em Mathematica}}

The recursive constructions underlying the proof of Butcher's
theorem can easily be realized using modern computer algebra
systems like {\em Mathematica}.\footnote{A {\em Mathematica}
notebook with all the programs and examples of this e-print is
included with the source files:
{http://arxiv.org/e-print/math.NA/0211049.tar.gz} .}
We assume that the reader is familiar with this particular
package.

We begin by defining the recursive data-structure of a rooted tree
as an unordered list---together with two simple routines for input and
output:
\smallskip\begin{flushleft}
{\tt Attributes[Tree]=\{Orderless,Listable\};\\
ToTree[f\_\,]:=Tree@@ToTree/@f; ToTree[f\_\,Symbol]:=Tree[];\\
Format[$\boldsymbol\beta$\_Tree]:=StringReplace[ToString[List@@$\boldsymbol\beta$/.\{\}->"$\boldsymbol\odot$"],\\
\qquad \{"\{"->"[","\}"->"]"\}]
}
\smallskip\end{flushleft}
\noindent
Now, here is an example for the input of
the tree representing the elementary differential $f''(f''(f,f'(f)),f) =
f''(f,f''(f'(f),f))$:
\smallskip\begin{flushleft}
{\tt $\boldsymbol\beta_{\tt 1}$=ToTree[f${^{\prime\prime}}$[f$^{\prime\prime}$[f,f$^{\prime}$[f]],f]];
$\boldsymbol\beta_{\tt
2}$=ToTree[f$^{\prime\prime}$[f,f$^{\prime\prime}$[f$^{\prime}$[f],f]]];
If[$\boldsymbol\beta_{\tt 1}$==$\boldsymbol\beta_{\tt
2}$,$\boldsymbol\beta_{\tt 1}$,,]}
\end{flushleft}
$$[\odot ,[\odot ,[\odot ]]]$$
\noindent The definition (\ref{eq:treedef}) of the order $\# \beta$ is simply
expressed
by the recursive procedure:

\smallskip\begin{flushleft}
{\tt TreeOrder[$\boldsymbol\beta\_\,$]:=1+Plus@@TreeOrder/@$\boldsymbol\beta$}
\smallskip\end{flushleft}

\noindent As an example, we take the elementary
differential $f''(f'''(f'(f),f'(f),f),f)$:
\smallskip
\begin{flushleft}
{\tt $\boldsymbol\beta$ = ToTree[f$^{\prime\prime}$[f$^{\prime\prime\prime}$[f$^{\prime}$[f],f$^{\prime}$[f],f],f]];
TreeOrder[$\boldsymbol\beta$]}
\end{flushleft}
$${8}$$

\noindent The definition (\ref{eq:factorialdef}) of $\beta !$ is analogously expressed
by:
\smallskip\begin{flushleft}
{\tt
 TreeFactorial[$\boldsymbol\beta$\_\,]:=TreeOrder[$\boldsymbol\beta$]Times@@TreeFactorial/@$\boldsymbol\beta$
}
\smallskip\end{flushleft}
The above used tree $\beta$ of order 8 gives:
 {\tt TreeFactorial[$\boldsymbol\beta$]
}
$${192}$$

\noindent For the sake of completeness we also express the recursive
definition~(\ref{eq:factorialdef}) of $\alpha_\beta$ using {\em Mathematica}:
\smallskip\begin{flushleft}
{\tt  TreeAlpha[$\boldsymbol\beta$\_\,]:=Length[Permutations[$\boldsymbol\beta$]]/Length[$\boldsymbol\beta$]!Times@@TreeAlpha/@$\boldsymbol\beta$}
\smallskip\end{flushleft}

\noindent An application to the above example yields: {\tt TreeAlpha[$\boldsymbol\beta$]}
$${\frac{1}{2}}$$

\noindent We are now ready to map the recursive definition (\ref{eq:abeta}) of
$\RKA^{(\beta)}$ and of the order condition $b^T \RKA^{(\beta)} = 1/\beta!$
to {\em Mathematica}:
\smallskip\begin{flushleft}
{\tt TreeA[$\boldsymbol\beta$\_,n\_\,:1]:=(vars=\{i,j,k,l,m,p,q,r,u,v,w\};\\
\quad Times@@(Sum[a$_{\tt{vars[[n]],vars[[n+1]]}}$TreeA[\#,n+1]//Evaluate, \\
\quad\qquad\qquad
\{vars[[n+1]],s\}//Evaluate]\&/@$\boldsymbol\beta$))\\
   TreeOrderCondition[$\boldsymbol\beta$\_\,]:=Sum[b$_{\tt i}$TreeA[$\boldsymbol\beta$]//Evaluate,\{i,s\}]\\
   \quad ==1/TreeFactorial[$\boldsymbol\beta$]
}
\smallskip\end{flushleft}

\noindent For convenience, the coordinate index of the vector
$\RKA^{(\beta)}$ can be chosen from the list
{\tt \{i,j,k,l,m,p,q,r,u,v,w\}} and is passed by number
as the second argument to {\tt TreeA}. The order condition
 belonging to the above example is obtained as follows:
 \smallskip\begin{flushleft}
 {\tt TreeOrderCondition[$\boldsymbol\beta$]
}
\end{flushleft}
$${\sum _{i=1}^{s}{b_i} \bigg(\sum _{j=1}^{s}{a_{i,j}}\bigg) \sum _{j=1}^{s}{a_{i,j}} \bigg(\sum _{k=1}^{s}{a_{j,k}}\bigg)
{{\Bigg(\sum _{k=1}^{s}{a_{j,k}} \sum
_{l=1}^{s}{a_{k,l}}\Bigg)}^2}==\frac{1}{192}}$$

\noindent Even the typesetting of this formula was done completely automatically, using
{\em Mathematica}'s ability to generate \TeX-sources.

To generate all the order conditions for a given order $p$, we need a device
that constructs the set of all trees $\beta$ with $\# \beta
\leq p$. There are, in principle, two different recursive approaches:
\begin{itemize}
\item {\em root-oriented}: generate all trees $\beta$ of order $\#
\beta = p$ by first, listing all integer partitions $p-1 = p_1 +
\ldots + p_n$, $n=1,\ldots,p-1$, and next, setting $\beta = [\beta_1,\ldots,\beta_n]$ for
all trees $\beta_1,\ldots,\beta_n$ of order $\# \beta_1 =
p_1 < p ,\ldots,\#\beta=p_n < p$. These trees have already been generated by
the recursion.
\item {\em leaf-oriented}: Add a leaf to each node of the
trees $\hat\beta = [\beta_1,\ldots,\beta_n]$ of
order $\# \hat\beta = p-1$, increasing thereby the order exactly by
one. This can be done recursively by adding a leaf to every
node of the subtrees $\beta_1,\ldots,\beta_n$.
\end{itemize}

\noindent The root-oriented approach was chosen by Sofroniou~\cite{Sofr94}
in his {\em Mathematica} package {\tt Butcher.m}. It requires an
efficient integer partition package and the handling of cartesian
products. The leaf-oriented approach is as least as efficient as the other one,
but much easier to code:

\smallskip\begin{flushleft}
{\tt
AddLeave[$\boldsymbol\beta$\_\,]:=AddLeave[$\boldsymbol\beta$]=Fold[Union,\{Prepend[$\boldsymbol\beta$,Tree[]]\},\\
  \quad
  ReplacePart[$\boldsymbol\beta$,AddLeave[$\boldsymbol\beta$[[\#]]],\#]\&/@Range[Length[$\boldsymbol\beta$]]]\\
 Trees[order\_\,]:=NestList[Union@@AddLeave/@\#\&,\{Tree[]\},order-1]
}\smallskip\end{flushleft}

\noindent Given an order $p$ this procedure generates a list of
the sets of trees for each order $q \leq p$, e.g., {\tt Trees[4]
}
\[
\{\{\odot \},\{[\odot ]\},\{[[\odot ]],[\odot , \odot ]\},
 \{[[[\odot ]]],[[\odot ,\odot ]],[\odot ,[\odot ]],[\odot ,\odot , \odot ]\}\}
\]
 \noindent For
instance, the number of trees for each order $p \leq 10$ is given
by the entries of the following list:
{\tt Length/@Trees[10] }
$${\{1,1,2,4,9,20,48,115,286,719\}}$$
 \noindent The number of
order conditions for
$p=10$ can thus be obtained by:
{\tt Plus@@\% }
$$1205$$

\noindent Finally, for concrete calculations one has to specify
the number $s$ of stages. The following procedure then generates the
specific set of equations for {\em explicit} Runge-Kutta methods:

\smallskip\begin{flushleft}
{\tt
 explicit=\{a$_{\tt i\_\,,j\_}$:>0/;i<=j,c$_{\tt 1}$->0\};\\
 OrderConditions[order\_,stages\_\,]:=(\\
\quad autonomous=Table[Sum[a$_{\tt i,j}$,\{j,stages\}]==c$_{\tt i}$,\{i,stages\}];\\
\quad \{(TreeOrderCondition/@Union@@Trees[order])/.s->stages/.\\
\quad ToRules[And@@autonomous],autonomous\}/.explicit)
}
\smallskip\end{flushleft}

\noindent This way, we can automatically generate and typeset the
order conditions for the classical explicit  4-stage  Runge-Kutta
methods of order~4:
\smallskip\begin{flushleft}
 {\tt First[OrderConditions[4,4]] }
\end{flushleft}
\[
\begin{split}
\big\{& {b_1}+{b_2}+{b_3}+{b_4}==1,\;{b_2}\multsp {c_2}+{b_3}
{c_3}+{b_4}\multsp {c_4}==\frac{1}{2},  \\*[2mm] & {b_3}\multsp
{c_2}\multsp {a_{3,2}}+{b_4}\multsp ({c_2}\multsp
{a_{4,2}}+{c_3}\multsp {a_{4,3}})==\frac{1}{6},\;{b_4}\multsp
{c_2}\multsp {a_{3,2}}\multsp {a_{4,3}}==\frac{1}{24},  \\*[2mm] &
{b_3}\multsp c_{2}^{2}\multsp {a_{3,2}}+{b_4}\multsp
(c_{2}^{2}\multsp {a_{4,2}}+c_{3}^{2}\multsp
{a_{4,3}})==\frac{1}{12},\;{b_2}\multsp c_{2}^{2}+{b_3}\multsp
c_{3}^{2}+{b_4}\multsp c_{4}^{2}==\frac{1}{3},  \\*[2mm] & {b_3}\multsp
{c_2}\multsp {c_3}\multsp {a_{3,2}}+{b_4}\multsp {c_4}\multsp
({c_2}\multsp {a_{4,2}}+{c_3}\multsp
{a_{4,3}})==\frac{1}{8},\;{b_2}\multsp c_{2}^{3}+{b_3}\multsp
c_{3}^{3}+{b_4}\multsp c_{4}^{3}==\frac{1}{4}\big\}
\end{split}
\]

\begin{Remark}
{\rm Even in the more recent literature one can find examples like
\cite{GaGu99},
where order conditions for Runge-Kutta methods are generated by using a
computer algebra system to calculate the Taylor
expansions of the flow and the discrete flow directly.
This approach is typically bound to {\em scalar}
non-autonomous equations, i.e., $d=1$. Besides being inefficient
for higher orders, it is well-known \cite{Butc63b} that for $p \geq 5$
 additional
order conditions for general systems make an appearance, which do
not show up in the scalar case.}
\end{Remark}

\vspace*{-0.5cm}
\section{Examples of usage}

The following simple procedure tempts to solve the order
conditions for a given order $p$ and stage number $s$ by using
brute force, i.e., {\em Mathematica}'s {\tt Solve}-com\-mand. To simplify the
task,
the user is allowed to supply a set {\tt pre} of a priori chosen
additional equations and assignments that he thinks to be helpful.

\smallskip\begin{flushleft}
{\tt  RungeKuttaMethod[p\_,s\_,pre\_\,]:=(\\
\quad rkTemplate=\{Table[a$_{\tt i, j}$,\{i,s\},\{j,s\}],Table[b$_{\tt
i}$,\{j,1\},\{i,s\}],\\
\quad\qquad Table[c$_{\tt i}$,\{j,1\},\{i,s\}]\}/.explicit;\\
\quad conditions=pre$\,\boldsymbol\cup\,$Union@@OrderConditions[p,s];\\
\quad solveVars=Complement[Flatten[rkTemplate],\{0\}];\\
\quad sol=Solve[conditions,solveVars];\\
\quad Thread[\{A,b,c\}==MatrixForm/@rkTemplate/.\#]\&/@sol)
}
\smallskip\end{flushleft}

\noindent Since {\em Mathematica}'s {\tt Solve}-command uses a Gr\"{o}bner-basis approach
for solving systems of polynomial equations, we can show this way that the classical
explicit 4-stage Runge-Kutta method of order $p=4$ is uniquely given by the
additional constraints $b_2=b_3$ and $c_2=c_3$:
{\tt RungeKuttaMethod[4,4,\{b$_{\tt
2}$==b$_{\tt 3}$,c$_{\tt
2}$==c$_{\tt 3}$\}] }
\[
\Big\{\Big\{A==\left(
        \begin{matrix}
        0&0&0&0 \\
        \frac{1}{2}&0&0&0 \\
        0&\frac{1}{2}&0&0 \\
        0&0&1&0
        \end{matrix} \right),b==\left(\begin{matrix}
        \frac{1}{6}&\frac{1}{3}&\frac{1}{3}&\frac{1}{6}
        \end{matrix}\right),c==\left(\begin{matrix}
        0&\frac{1}{2}&\frac{1}{2}&1
  \end{matrix} \right)\Big\}\Big\}
\]

\noindent The next example is more demanding. In his book \cite[p.
199]{Butc87}, Butcher describes an algorithm for the construction
of explicit 6-stage methods of order $p=5$. The choices $c_6=1$
and $b_2=0$ together with the free parameters $c_2,c_3,c_4,c_5$ and
$a_{43}$ yield a unique method. Butcher provides a two-parameter
example by choosing $c_2 = u,c_3=1/4,c_4=1/2,c_5=3/4,a_{43}=v$. By
just passing this additional information to {\em Mathematica}'s {\tt
Solve}-command we obtain the following solution:
\smallskip\begin{flushleft}{\tt
 pre = \{c$_{\tt 2}$==u,c$_{\tt 3}$==1/4,c$_{\tt 4}$==1/2,c$_{\tt 5}$==3/4,c$_{\tt 6}$==1,b$_{\tt
2}$==0,a$_{\tt 4,3}$==v\};\\
 First[RungeKuttaMethod[5,6,pre]] }
\end{flushleft}
\[
\begin{split}
\Big\{& A==\left(\begin{matrix}
    0&0&0&0&0&0 \\*[2mm]
    u&0&0&0&0&0 \\*[2mm]
    \frac{-1+8\multsp u}{32\multsp u}&\frac{1}{32\multsp u}&0&0&0&0 \\*[2mm]
    \frac{-1+4\multsp u+2\multsp v-8\multsp u\multsp v}{8\multsp u}&
    \frac{1-2\multsp v}{8\multsp u}&v&0&0&0 \\*[2mm]
    \frac{3\multsp (1-3\multsp u-v+4\multsp u\multsp v)}{16\multsp u}&
    \frac{3\multsp (-1+v)}{16\multsp u}&
    -\frac{3}{4}\multsp (-1+v)&\frac{9}{16}&0&0\\*[2mm]
    \frac{-7+22\multsp u+6\multsp v-24\multsp u\multsp v}{14\multsp u}&
    \frac{7-6\multsp v}{14\multsp u}&\frac{12\multsp v}{7}&-\frac{12}{7}&\frac{8}{7}&0
  \end{matrix} \right)\\*[2mm]
& b==\left(\begin{matrix}
    \frac{7}{90} &
    0 &
    \frac{16}{45} &
    \frac{2}{15} &
    \frac{16}{45} &
    \frac{7}{90}
   \end{matrix} \right),\quad c==\left(\begin{matrix}
    0 &
    u &
    \frac{1}{4} &
    \frac{1}{2} &
    \frac{3}{4} &
    1
  \end{matrix} \right)\quad \Big\}
  \end{split}
\]

\noindent This result shows that the coefficients $a_{51}$ and $a_{52}$ of
Butcher's solution \cite[p.
199]{Butc87} are in error, a fact that was already observed
by Sofroniou~\cite{Sofr94} using the {\em Mathematica} package {\tt Butcher.m}.

\providecommand{\bysame}{\leavevmode\hbox to3em{\hrulefill}\thinspace}
\providecommand{\MR}{\relax\ifhmode\unskip\space\fi MR }
\providecommand{\MRhref}[2]{%
  \href{http://www.ams.org/mathscinet-getitem?mr=#1}{#2}
}
\providecommand{\href}[2]{#2}


\begin{thebibliography}{1}

\bibitem{Butc63}
J.~C. Butcher, \emph{Coefficients for the study of {Runge-Kutta} integration
  processes}, J. Austral. Math. Soc. \textbf{3} (1963), 185--201.

\bibitem{Butc63b}
\bysame, \emph{On the integration processes of {A. Huta}}, J. Austral. Math.
  Soc. \textbf{3} (1963), 202--206.

\bibitem{Butc87}
\bysame, \emph{The numerical analysis of ordinary differential equations.
  {Runge-Kutta} methods and general linear methods}, {John Wiley \& Sons},
  Chichester, 1987.

\bibitem{GaGu99}
W.~Gander and D.~Gruntz, \emph{Derivation of numerical methods using computer
  algebra}, SIAM Review \textbf{41} (1999), 577--593.

\bibitem{Hair99}
E.~Hairer, \emph{Numerical geometric integration. {Lecture notes of a course
  given in 1998/99}}, {\tt http://www.unige.ch/math/folks/hairer/polycop.html},
  1999.

\bibitem{HNWa87}
E.~Hairer, S.~P. {N\o rsett}, and G.~Wanner, \emph{Solving ordinary
  differential equations i. nonstiff problems}, 2. ed., {Springer-Verlag},
  Berlin, Heidelberg, New York, 1993.

\bibitem{HaWa74}
E.~Hairer and G.~Wanner, \emph{On the {Butcher} group and general multi-value
  methods}, Computing \textbf{13} (1974), 1--15.

\bibitem{Jenk76}
R.~J. Jenks, \emph{Problem \# 11: Generation of {Runge-Kutta} equations},
  SIGSAM Bulletin \textbf{10} (1976), 6.

\bibitem{Sofr94}
M.~Sofroniou, \emph{Sysmbolic derivation of {Runge-Kutta} methods}, J. Symb.
  Comp. \textbf{18} (1994), 265--296.

\end{thebibliography}
\end{document}